\documentclass[12pt]{article}

\usepackage{graphicx} 
\usepackage{amsmath} 
\usepackage{amssymb,amsthm}
\usepackage{color}
\usepackage{tikz}
\usepackage{pgfplots}

\newcommand{\R}{{\mathbb R}}\newcommand{\N}{{\mathbb N}}

\newcommand{\C}{{\mathbb C}}

\let\epsilon\varepsilon

\let\theta\vartheta

\newtheorem{theorem}{Theorem}[section]\newtheorem{lemma}[theorem]{Lemma}

\newtheorem{remark}[theorem]{Remark}

\title{A unified Sobolev framework for long-wave limits of the Gross-Pitaevskii equation}
\author{Guido Schneider, Nils Thorin\footnote{Corresponding author: \texttt{nils.thorin@mathematik.uni-stuttgart.de}} \\
{\small 
Institut f\"ur Analysis, Dynamik und Modellierung,} \\ 
{\small Universit\"at Stuttgart, Pfaffenwaldring 57, 70569 Stuttgart, Germany} \\
}

\begin{document}


\maketitle
\begin{abstract}
We present a unified justification of the KdV, Airy, inviscid Burgers, and Whitham approximations for the one-dimensional Gross-Pita\-evskii equation in Sobolev spaces. Using energy estimates for its Madelung formulation, we establish quantitative error bounds on the natural time scales of the respective approximations. 
The analysis covers different long-wave scaling regimes within a single framework.
\end{abstract}

\section{Introduction}

We consider the one-dimensional Gross-Pitaevskii (GP) equation
\begin{equation} \label{GP}
i \partial_t q = - \partial_x^2 q +2 (|q|^2-1)q,
\end{equation}
with $x,t \in \R $ and $ q(x,t) \in \C$.
We are interested in the description of the dynamics near $ q = 1$ 
by amplitude equations in various  long wave limits.
For sufficiently small amplitudes, linear dispersion dominates and leads to the Airy equation. At the critical scaling, dispersion and nonlinear transport balance, yielding the KdV equation. For larger but still asymptotically vanishing amplitudes, nonlinear transport dominates, resulting in the inviscid Burgers equation. When the amplitude remains fixed in the long-wave limit, the coupled evolution is described by Whitham’s modulation equations. We prove the validity of these approximations within a unified Sobolev framework.

In this introduction we first introduce the Madelung formulation of the Gross-Pitaevskii equation and derive the amplitude equations arising in the different long-wave regimes. We then state the approximation results in Sobolev spaces, discuss their relation to previous work, and outline the common proof strategy and the organization of the paper.

\subsection{The Madelung transform}

For the derivation of these amplitude equations it is
advantageous to apply 
the Madelung transform $ q \mapsto (\rho,v) $. Assuming $\rho=|q|^2>0$, we write
\begin{equation}
q = \sqrt{\rho} e^{i \phi}, \qquad v = \partial_x \phi.
\end{equation}
Using 
\begin{eqnarray*}
i \partial_t q & = & i \frac12 \frac{\partial_t \rho}{\sqrt{\rho}} e^{i \phi} -  \sqrt{\rho} e^{i \phi} \partial_t \phi , \\
- \partial_x^2 q & = & - \frac12 \frac{\partial_x^2 \rho}{\sqrt{\rho}} e^{i \phi} 
 + \frac14 \frac{(\partial_x \rho)^2}{\rho^{3/2}} e^{i \phi} 
- i  \frac{\partial_x \rho}{\sqrt{\rho}} e^{i \phi} \partial_x \phi
+  \sqrt{\rho} e^{i \phi} (\partial_x \phi)^2
-i  \sqrt{\rho} e^{i \phi} \partial_x^2 \phi,
\end{eqnarray*}
and
$$
2 (|q|^2-1)q  =  2 (\rho- 1) \sqrt{\rho} e^{i \phi},
$$
and separating the real and imaginary parts, then differentiating the phase
equation with respect to $x$, we obtain
\begin{eqnarray*}
\partial_t \rho + 2 \partial_x (\rho v) & = & 0 , \\
\partial_t v + \partial_x (v^2)+ 2 \partial_x \rho & = & \partial_x(\partial_x(\frac{1}{2} \frac{\partial_x \rho}{\rho})+ (\frac{1}{2} \frac{\partial_x \rho}{\rho})^2 ).
\end{eqnarray*}
Using
$$
 \frac{\partial_x \rho}{\rho} = \partial_x (\ln \rho),
$$
introducing $ w = \ln \rho $,  and
using
$
\partial_t w = \frac{\partial_t \rho}{\rho}
$,
gives  finally 
\begin{eqnarray}\label{made1}
\partial_t w & = & - 2 \partial_x v - 2 v \partial_x w  , \\
\partial_t v & = &  \frac12 \partial_x^3 w - \partial_x (v^2)- 2 \partial_x (e^w) +
 \frac14 \partial_x ((\partial_x w)^2).
\label{made2}
\end{eqnarray}
We are now interested in the dynamics near the point $ (v,w) = (0,0) $. In particular, we are interested in various long-wave limits of this system.

\subsection{The long wave approximations}

\label{sec12}

The amplitude equations will be derived via multiple scaling analysis
using a small perturbation parameter $ 0 < \delta \ll 1 $. 
We will derive  the Airy equation, the KdV equation, the inviscid Burgers equation, and Whitham’s modulation equations.

i) For the derivation of the KdV equation we make the ansatz 
$$
v(x,t) = \delta^2 V(\delta (x-ct),\delta^3 t) , \qquad 
w(x,t) = \delta^2 W(\delta (x-ct),\delta^3 t).
$$
We find
\begin{eqnarray*}
\delta^3 (-c \partial_X W) + \delta^5 \partial_T W & = & - 2 \delta^3 \partial_X V - 2 \delta^5 V \partial_X W  , \\
\delta^3 (-c \partial_X V) + \delta^5 \partial_T V & = &  \frac12 \delta^5 \partial_X^3 W - \delta^5 \partial_X (V^2)- 2 \delta \partial_X (e^{\delta^2 W})  + \mathcal{O}(\delta^7).
\end{eqnarray*}
Using 
$$ 
- 2 \delta \partial_X (e^{\delta^2 W}) = - 2 \delta^3 \partial_X W -  \delta^5 \partial_X (W^2) 
+ \mathcal{O}(\delta^7)
$$ 
and equating the coefficients in front of $ \delta^3 $
to zero gives 
$$ 
-c \partial_X W = -2  \partial_X V , \qquad  -c \partial_X V = - 2  \partial_X W.
$$ 
This yields $c=\pm2$ and
$V=\frac{c}{2}W+C(T)$. Assuming that $ V $ and $ W $ decay to zero as
$|X|\to\infty$, the integration constant $ C(T) $ vanishes, giving the leading-order relations
$$ 
c = 2 , \quad V = W \qquad \textrm{or} \qquad c = - 2 , \quad V = - W.
$$ 
Equating the coefficients in front of $ \delta^5 $
to zero we obtain 
in case $ 
c = 2 $ and  with $  V = W + a \delta^2 \partial_X^2 W + b \delta^2 W^2 $ that 
\begin{eqnarray*}
   \partial_T W & = &  
- 2 a \partial_X^3 W -2  b \partial_X (W^2) 
- 2  W \partial_X W  , \\
- 2 a \partial_X^3 W -2  b  \partial_X (W^2) + \partial_T W & = &  \frac12  \partial_X^3 W -  \partial_X (W^2)-  \partial_X (W^2) .
\end{eqnarray*}
Subtracting the first equation from the second equation gives for the coefficient 
of $ \partial_X^3 W $ the condition
$ -2  a = 1/2 + 2 a $ and for the coefficient 
of $ \partial_X (W^2) $ the condition $ -2 b = -2 + 2 b + 1$. 
Therefore, we obtain $ a = -1/8 $ and $ b = 1/4 $ and so
we finally find the KdV equation 
\begin{equation} \label{KdV}
\partial_T W  =   \frac14  \partial_X^3 W  -\frac32  \partial_X (W^2) .
\end{equation}

ii)+iii) This ansatz can be generalized to 
$$
v(x,t) = \delta^{\alpha} V(\delta (x-ct),\delta^{1+\alpha} t) , \qquad 
w(x,t) = \delta^{\alpha} W(\delta (x-ct),\delta^{1+\alpha} t)
$$
with $ \alpha > 0  $. 
We find 
\begin{eqnarray}\label{made1a}
-c \delta^{1+\alpha} \partial_X W+ \delta^{1+2\alpha}\partial_T W & = & - 2 \delta^{1+\alpha} \partial_X V - 2 \delta^{1+ 2 \alpha}  V \partial_X W  , \\
\label{made2a}
-c \delta^{1+\alpha} \partial_X V+ \delta^{1+2\alpha}\partial_T V & = & \delta^{3+\alpha}  \frac12 \partial_X^3 W - \delta^{1+2 \alpha}  \partial_X (V^2)
\\ &&- 2 \delta \partial_X (e^{\delta^{\alpha}W}) +
 \frac14 \delta^{3+ 2 \alpha} \partial_X ((\partial_X W)^2).
\nonumber
\end{eqnarray}
As above,  we obtain the leading-order relations
$$ 
c = 2 , \quad V = W \qquad \textrm{or} \qquad c = - 2 , \quad V = - W.
$$ 
For $\alpha>2$, dispersion enters at order $\delta^{3+\alpha}$,
before the time derivative at order $\delta^{1+2\alpha}$ for the above time scale.
We therefore replace the slow time by $T=\delta^3t$ and use
$$
v(x,t)=\delta^\alpha V(\delta(x-2t),\delta^3t), \qquad
w(x,t)=\delta^\alpha W(\delta(x-2t),\delta^3t).
$$
With $V=W+a\delta^2\partial_X^2W$, matching terms at order
$\delta^{3+\alpha}$ gives $a=-1/8$ and the Airy equation
\begin{equation}\label{Airy}
\partial_T W = \frac14  \partial_X^3 W 
\end{equation}
at leading order. On this time scale, the nonlinear terms are smaller by
a factor $\delta^{\alpha-2}$.
The KdV equation is recovered for $ \alpha = 2 $. For $ \alpha \in (0,2) $,
we keep the original slow time $T=\delta^{1+\alpha}t$ and
we choose  $ c = 2 $ and $ V = W + b \delta^{\alpha} W^2$. We find again $ b= 1/4 $
and therefore  the inviscid Burgers equation 
\begin{equation}
\label{constburgers1}
\partial_T W =  -\frac32  \partial_X (W^2) 
\end{equation}
is obtained.

iv) By making the ansatz 
$$
v(x,t) = V(\delta x,\delta t) , \qquad 
w(x,t) = W(\delta x,\delta t),
$$
we directly obtain Whitham's modulation equations 
\begin{eqnarray}\label{wh1}
\partial_T W & = & - 2 \partial_X V - 2 V \partial_X W  , \\
\partial_T V & = &  - \partial_X (V^2)-  2 \partial_X (e^W) ,
\label{wh2}
\end{eqnarray}
by equating the coefficients in front of $ \delta^1 $
to zero.

%
%

\subsection{The approximation results}

Our approximation results  measure the errors in Sobolev norms on the natural time scales of the  amplitude equations. 

\begin{remark}{\rm 
Rigorous KdV approximation results for the Gross-Pitaevskii equation were established in \cite{BGSS09,BGSS10,CR10}. These works employ integrability properties, energy estimates, and hydrodynamic formulations. A Cauchy-Kowalevskaya-type approach in Gevrey  spaces 
using the 
conservation law structure of the problem was used in \cite{CDS14}, while \cite{We26} treats long-wave approximation at the level of the NLS and KdV hierarchies. Whitham's modulation approximation results were obtained in \cite{DS09,BKZ21}. Our method is particularly close to \cite{BKZ21}, where higher-order corrections and energy estimates yield validity in Sobolev spaces. A common treatment of the KdV, Burgers, and Whitham's modulation limits was developed in \cite{BDS19} for a spatially periodic Boussinesq model. Here, we establish a unified Sobolev framework for the KdV, Airy, inviscid Burgers, and Whitham's modulation  approximations of the GP equation. To our knowledge, the Airy and inviscid Burgers justification results are not documented  in the existing literature.
}\end{remark}
For  the formulation of our approximation results we 
use classical Sobolev spaces $ H^m $ for $ m \in \N_0 $.
They are equipped with the norm
$$
\| u \|_{H^m} = \sum_{j=0}^m \| \partial_x^j u \|_{L^2}, \qquad \textrm{with} 
\qquad \| u \|^2_{L^2} = \int_{-\infty}^{\infty} |u(x)|^2 dx. 
$$
Our approximation results for the KdV, Airy, inviscid Burgers, and Whitham's modulation approximations
are as follows:
\begin{theorem} \label{th11}
For all $ m \in \N_0 $ with $ m \geq 2 $ the following holds.
Let  $W \in C([0,T_0],H^{m+5})$ be a solution of the KdV equation \eqref{KdV}. Then there exist $ \delta_0 > 0 $ and $ C > 0 $ such that
for all $ \delta \in (0,\delta_0)$ we have solutions $ (v,w) $ of \eqref{made1}-\eqref{made2} with 
$$ 
\sup_{t \in [0,T_0/\delta^3]} \| (v,w)(\cdot,t) - \delta^2 W(\delta(\cdot-2 t), \delta^3 t )(1,1) \|_{H^m \times H^{m+1}} \leq C \delta^{7/2}. 
$$
\end{theorem}
\begin{remark}{\rm
The power in the error bound reflects the first neglected terms in the KdV
approximation and the scaling of the Sobolev norm.
With $V=W-\frac18\delta^2\partial_X^2W+\frac14\delta^2W^2$,
the formal residual is of order $\delta^7$.
Over times of order $\delta^{-3}$ this gives an error of order $\delta^4$.
Changing to the unscaled Sobolev norm gives a factor $\delta^{-1/2}$,
leading to the bound $C\delta^{7/2}$.
}\end{remark}


For amplitudes smaller than those in the KdV scaling, we obtain the Airy approximation.
The pure Airy approximation from above is not sufficient to formulate a satisfying 
approximation result. For the formulation and the proof the above 
Airy approximation has to be extended by higher order terms.
The computation of the higher approximation terms is straightforward.
They can be computed as long as the solution of the Airy equation exists.
\begin{theorem}\label{thairy}
For all $\alpha>2$ and $m \in\N_0$ with $ m \geq 2 $ the following holds.
Let $T_0>0$ and let $W\in C([0,T_0],H^{m_W})$ be a solution of the Airy equation
\eqref{Airy}, with $ m_W \geq m $ sufficiently large. Then there exist $\delta_0>0$ and $C>0$
such that for all $\delta\in(0,\delta_0)$ there is a higher order approximation
$$ 
\delta^\alpha (V_h,W_h)(X,T) = \delta^\alpha W(X,T)(1,1) 
+\mathcal O\!\left(\delta^{\min\{2\alpha-2,\alpha+2\}}\right)
$$
for which there exist solutions $(v,w)$ of \eqref{made1}-\eqref{made2} with
$$
\sup_{t\in[0,T_0/\delta^3]}
\left\|(v,w)(\cdot,t)- \delta^\alpha (V_h,W_h)(\delta(\cdot-2t),\delta^3t)
\right\|_{H^m\times H^{m+1}}
\leq C\delta^{\alpha+ 3/2}.
$$
\end{theorem}

The same problem appears for the inviscid Burgers approximation from above.
 The pure inviscid Burgers approximation is again  not sufficient to formulate a satisfying 
approximation result. For the formulation and the proof the above 
inviscid Burgers approximation has to be extended by higher order terms.
The computation of the higher approximation terms is straightforward.
They can be computed as long as the solution of the inviscid Burgers  equation exists.
\begin{theorem} \label{th12}
For all $ \alpha \in (0,2) $ and $ m \in \N_0 $ with $ m \geq 2 $ the following holds.
Let   $W \in C([0,T_0],H^{m_W})$, with $ m_W \geq m $ sufficiently large, be a solution of the inviscid Burgers  equation \eqref{constburgers1}. Then there exist $ \delta_0 > 0 $ and $ C > 0 $ such that
for all $ \delta \in (0,\delta_0) $ 
there is a higher order approximation
$$ 
\delta^\alpha (V_h,W_h)(X,T) = \delta^\alpha W(X,T)(1,1) 
+\mathcal O\!\left(\delta^{\min\{2\alpha,2\}}\right)
$$
for which
there exist solutions $ (v,w) $ of \eqref{made1}-\eqref{made2}  with 
$$ 
\sup_{t \in [0,T_0/\delta^{1+\alpha}]} \| (v,w)(\cdot,t) - \delta^\alpha (V_h,W_h)(\delta(\cdot-2 t), \delta^{1+\alpha} t ) \|_{H^m \times H^{m+1}} \leq C \delta^{\alpha+3/2}. 
$$
\end{theorem}

Finally,  Whitham’s modulation equations describe modulations 
of order $ \mathcal{O}(1)$. 
\begin{theorem} \label{th11whit}
For all $ m \in \N_0 $ with $ m \geq 2 $
and all $ C_1 > 0 $  the 
following holds.  Let $(V,W) \in C([0,T_0],H^{m+3}\times H^{m+3})$ be a solution of   \eqref{wh1}-\eqref{wh2}
with 
\begin{equation} \label{condition1}
 \sup_{T \in [0,T_0]} \| (V,W)(\cdot,T) \|_{H^{m+3}\times H^{m+3}} \leq C_1 .
 \end{equation}
 Then there exist $ \delta_0 > 0 $ and $ C > 0 $, depending on $m$, $T_0$, and $C_1$, but independent of $\delta$, such that
for all $ \delta \in (0,\delta_0)$ we have solutions $ (v,w) $ of \eqref{made1}-\eqref{made2}  with 
$$ 
\sup_{t \in [0,T_0/\delta]} \| (v,w)(\cdot,t) - (V,W)(\delta \cdot, \delta t ) \|_{H^m \times H^{m+1}} \leq C \delta^{3/2}. 
$$
\end{theorem}
\begin{remark}{\rm
The approximation result on Whitham’s modulation equations, stated in \cite{BDS19}, contains 
a smallness condition in the sense that the formulation 'for all $ C_1 > 0 $'
is replaced by 'there exists a $ C_1 > 0 $'. 
For the Gross-Pitaevskii equation
such a condition is not necessary due to the special structure of the problem.
See Remark \ref{remequi}.
This is the reason why we have chosen a slightly different formulation of 
Theorem \ref{th11whit} compared to the theorems above.
}\end{remark}
\begin{remark}{\rm
Whitham’s modulation approximation keeps the full nonlinear terms and neglects the
terms containing three spatial derivatives such that we have a formal residual of order $\delta^3$
in the original equations. Over times of order $\delta^{-1}$, formal
integration therefore suggests an error of order $\delta^2$, or
$\delta^{3/2}$ in the unscaled Sobolev norm.
This is the rate stated in the theorem. 
}\end{remark}

\begin{remark}{\rm 
Proving that the formal KdV, Airy, Burgers, or Whitham’s modulation  approximations remain close to true solutions of the original system over their natural time scales
is a non-trivial task, since solutions with amplitude $\mathcal{O}(\delta^{\alpha})$ must be
controlled on time intervals of length $\mathcal{O}(\delta^{-\min(1+\alpha,3)})$. 
}\end{remark}

\begin{remark}{\rm 
A discussion about the transfer of our approximation results
to the  variable $ q $ of  the original GP equation \eqref{GP}
can be found in Section 
\ref{secdisc}.
}
\end{remark}

%

The plan of the paper is as follows. In Section \ref{sec2}
we construct a higher order approximation in order to make the residual,
the terms which do not cancel after inserting the approximation, 
as small as wanted. The energy estimates 
for the error are proven in Section \ref{sec3}. The paper is closed with an outlook and a discussion.
\medskip

{\bf Notation.}
Many different constants are denoted with the same symbol $ C $ if they can be chosen independently of the small perturbation parameter $ 0 < \delta \ll 1 $.
Unless indicated otherwise, we will use the abbreviation $ \int = \int_{-\infty}^{\infty} $.

\section{Estimates for the residual terms}
\label{sec2}

The residual
\begin{eqnarray}\label{res1}
 Res_w & = & -\partial_t w-2\partial_x v-2v\partial_x w, \\
 Res_v & = & -\partial_t v+\frac12\partial_x^3 w-\partial_x(v^2)
 -2\partial_x(e^w)+\frac14\partial_x((\partial_x w)^2)
\label{res2}
\end{eqnarray}
measures how much the KdV, Airy, inviscid Burgers, and Whitham’s modulation  approximations 
fail to satisfy \eqref{made1}-\eqref{made2}. The residual appears as a forcing term
in the error equations. All residuals below are defined in the original
variables $(x,t)$ and are measured in the unscaled Sobolev norms.
For $0<\delta\leq1$, we use the scaling identity
\begin{equation}\label{rescalingHs}
\|\delta^p F(\delta(\cdot-ct))\|_{H^r}
=\sum_{j=0}^r\delta^{p+j-1/2}\|\partial_X^jF\|_{L^2}
\leq\delta^{p-1/2}\|F\|_{H^r},\qquad r\in\N_0.
\end{equation}

i) We start with the KdV approximation.
\begin{lemma}\label{lemreskdv}
Let $m\in\N_0$, $T_0>0$, and let $W\in C([0,T_0],H^{m+5})$
solve \eqref{KdV}. Set $X=\delta(x-2t)$, $T=\delta^3t$, and
$$
w(x,t)=\delta^2W(X,T),\qquad
v(x,t)=\delta^2\left(W-\frac18\delta^2\partial_X^2W
+\frac14\delta^2W^2\right)(X,T).
$$
Then there exist $\delta_0\in(0,1]$ and $C>0$ such that
for all $\delta\in(0,\delta_0)$ we have
$$
\sup_{t\in[0,T_0/\delta^3]}
\|(Res_v,Res_w)(\cdot,t)\|_{H^m\times H^{m+1}}
\leq C\delta^{13/2}.
$$
\end{lemma}
\begin{proof}
Write $V=W+a\delta^2\partial_X^2W+b\delta^2W^2$, with $a=-1/8$
and $b=1/4$. Substitution and \eqref{KdV} cancel the terms of orders
$\delta^3$ and $\delta^5$ as seen in the introduction. The remaining terms start at order $\delta^7$.
The term with the highest spatial derivatives in $Res_v$ comes from $-\partial_t v$:
$$
-a\delta^7\partial_X^2\partial_TW
=-a\delta^7\partial_X^2\left(\frac14\partial_X^3W
-\frac32\partial_X(W^2)\right).
$$
It is bounded in $H^m$ by the assumed $H^{m+5}$ bound on $W$.
The remaining terms in $Res_v$ and $Res_w$ are uniformly bounded
in $H^m$ and $H^{m+1}$, respectively, by Sobolev product estimates
and Taylor's formula for the exponential. 
The factor $\delta^{-1/2}$ in these residual bounds comes from the
$L^2$ scaling in \eqref{rescalingHs}.
\end{proof}

ii) The Airy approximation 
from the introduction yields a residual which is too large for our purposes.
Therefore, we construct a higher order Airy approximation with the following properties. 
\begin{lemma}\label{lemresairy}
Let $\alpha>2$, $m\in\N_0$, $T_0>0$, and let
$W\in C([0,T_0],H^{m_W})$ solve \eqref{Airy}, with $ m_W \geq  m  $ sufficiently large.
Set $X=\delta(x-2t)$ and  $T=\delta^3t$. Then there exists a higher order approximation 
\begin{eqnarray*}
w(x,t)& = & \delta^\alpha W_h(X,T) = \delta^\alpha W(X,T)
+\mathcal O\!\left(\delta^{\min\{2\alpha-2,\alpha+2\}}\right)
,\\
v(x,t)& =& \delta^\alpha\ V_h(X,T) = \delta^\alpha W(X,T)
+\mathcal O\!\left(\delta^{\min\{2\alpha-2,\alpha+2\}}\right),
\end{eqnarray*}
for which there exist $\delta_0\in(0,1]$ and $C>0$ such that
for all $\delta\in(0,\delta_0)$ we have
$$
\sup_{t\in[0,T_0/\delta^3]}
\|(Res_v,Res_w)(\cdot,t)\|_{H^m\times H^{m+1}}
\leq C\delta^{\alpha + 9/2}.
$$
\end{lemma}
\begin{proof}
The proof is straightforward, but the number of equations depends on $ \alpha $ and increases for $ \alpha \to 2 $. 
Since the construction of the higher order approximation is not the main purpose of this paper, we restrict our presentation of this topic  to the case $ \alpha = 3 $ for which we exemplarily 
explain the procedure in Appendix \ref{appHigherAiry}. 
\end{proof}

iii) The  inviscid Burgers approximation
from the introduction yields a residual which is too large for our purposes.
Therefore, we construct a higher order inviscid Burgers  approximation with the following properties. 
\begin{lemma}\label{lemresburgers}
Let $\alpha\in(0,2)$, $m\in\N_0$, $T_0>0$, and let
$W\in C([0,T_0],H^{m_W})$, with $ m_W \geq  m  $ sufficiently large,  solve \eqref{constburgers1}.
Set $X=\delta(x-2t)$ and  $T=\delta^{1+\alpha}t$. Then there exists a higher order approximation 
\begin{eqnarray*}
w(x,t)& = & \delta^\alpha W_h(X,T) = \delta^\alpha W(X,T)
+\mathcal O\!\left(\delta^{\min\{2\alpha,2\}}\right),\\
v(x,t)& =& \delta^\alpha V_h(X,T) = \delta^\alpha W(X,T)
+\mathcal O\!\left(\delta^{\min\{2\alpha,2\}}\right),
\end{eqnarray*}
for which there exist $\delta_0\in(0,1]$ and $C>0$ such that
for all $\delta\in(0,\delta_0)$ we have 
$$
\sup_{t\in[0,T_0/\delta^{1+\alpha}]}
\|(Res_v,Res_w)(\cdot,t)\|_{H^m\times H^{m+1}}
\leq C\delta^{2 \alpha+5/2}.
$$
\end{lemma}
\begin{proof}
The proof is straightforward, but the number of equations depends on $ \alpha $ and 
increases as $ \alpha \to 0 $ or $ \alpha \to 2 $. 
Since the construction of the higher order approximation is not the main purpose of this paper, we restrict our presentation of this topic  to the case $ \alpha = 1 $ for which we exemplarily 
explain the procedure in Appendix \ref{appHigherBurgers}. 
\end{proof}

iv) Finally we consider Whitham’s modulation approximation.
\begin{lemma}\label{lemreswme}
Let $m\in\N_0$, $T_0>0$, and let
$(V,W)\in C([0,T_0],H^{m+3}\times H^{m+3})$ solve \eqref{wh1}-\eqref{wh2}, with
$$
\sup_{T\in[0,T_0]}\|(V,W)(\cdot,T)\|_{H^{m+3}\times H^{m+3}}\leq C_1.
$$
Set $v(x,t)=V(\delta x,\delta t)$ and $w(x,t)=W(\delta x,\delta t)$.
There exists $C>0$, depending only on $m$ and $C_1$, such that
for all $\delta\in(0,1]$ we have 
$$
\sup_{t\in[0,T_0/\delta]}
\|(Res_v,Res_w)(\cdot,t)\|_{H^m\times H^{m+1}}
\leq C\delta^{5/2}.
$$
\end{lemma}
\begin{proof}
Whitham’s modulation equations give the exact identities
\begin{eqnarray*}
Res_w(x,t)&=&0,\\
Res_v(x,t)&=&\delta^3\left(\frac12\partial_X^3W
+\frac14\partial_X((\partial_XW)^2)\right)(\delta x,\delta t).
\end{eqnarray*}
The term in parentheses is uniformly bounded in $H^m$ by the assumed
$H^{m+3}$ bound and Sobolev product estimates.
The factor $\delta^{-1/2}$ in these residual bounds comes from the
$L^2$-scaling in \eqref{rescalingHs}.
\end{proof}

\section{The error estimates}

\label{sec3}

We now turn to the main objective of this paper, namely a unified Sobolev framework that enables us to derive error bounds for the various long-wave approximations constructed in the previous section simultaneously for all values of $ \alpha \geq 0 $.

Our starting point is the system \eqref{made1}-\eqref{made2} which we recall
\begin{eqnarray}\label{madeErrorW}
\partial_t w & = & - 2 \partial_x v - 2 v \partial_x w  , \\
\partial_t v & = &  \frac12 \partial_x^3 w - \partial_x (v^2)- 2 \partial_x (e^w) +
 \frac14 \partial_x ((\partial_x w)^2).
\label{madeErrorV}
\end{eqnarray}
The errors  $\delta^\beta R_V$ and $\delta^\beta R_W$ made by 
the approximations  $ \delta^\alpha V(X,T) $ and $ \delta^\alpha W(X,T) $
are defined by
\begin{eqnarray} \label{errordef1}
v(x,t)&=&\delta^\alpha V(X,T)+\delta^\beta R_V(x,t),\\
w(x,t)&=&\delta^\alpha W(X,T)+\delta^\beta R_W(x,t). \label{errordef2}
\end{eqnarray}
Here, with a slight abuse of notation, $V$ and $W$ denote the full approximations, including their
$\delta$-dependent corrections. 
Moreover, we have  $X=\delta(x-ct)$ and $T$, $ c $, and $ \beta $ 
according to the following table:
\[
\begin{array}{|c|c|c|c|c|} \hline
\text{Approximation} & \alpha & c & T & \beta \\ \hline
\text{KdV} & 2 & 2 & \delta^3t & \frac72 \\
\text{Higher order Airy} & \alpha>2 & 2 & \delta^3t
 & \alpha+\frac32 \\
\text{Higher order Burgers} & 0<\alpha<2 & 2 & \delta^{1+\alpha}t
 & \alpha+\frac32 \\
\text{Whitham’s modulation} & 0 & 0 & \delta t & \frac32 \\ \hline
\end{array}
\]
The desired error estimates correspond to bounds on
$(R_V,R_W)$ in $H^m\times H^{m+1}$ that are uniform in $\delta$.

In the following equations, $V$, $W$, and $\partial_XW$ are evaluated at
$(X,T)$, while $R_V$ and $R_W$ are evaluated at $(x,t)$.
The derivative $\partial_x$ acts on the composed functions; in particular,
$\partial_x W(X,T)=\delta\partial_XW(X,T)$.
The residuals are those of the corresponding corrected approximations.

After inserting \eqref{errordef1}-\eqref{errordef2} into 
\eqref{madeErrorW}-\eqref{madeErrorV} and dividing everything by $\delta^\beta$,
we find
\begin{eqnarray*}
\partial_t R_W & = & - 2 \partial_x R_V 
 - 2 \delta^{\alpha}  V \partial_x R_W
 - 2 \delta^{1+  \alpha }  R_V \partial_X W
\\ && - 2 \delta^{\beta}  R_V \partial_x R_W
  +  \delta^{-\beta}Res_w , \\
\partial_t R_V & = &  \frac12 \partial_x^3 R_W - 2 \delta^{\alpha}  \partial_x (V R_V) - \delta^{\beta}  \partial_x (R_V^2)
\\ && +
 \frac12 \delta^{1+ \alpha} \partial_x ((\partial_X W)(\partial_x R_W))
 +
 \frac14 \delta^{\beta} \partial_x ((\partial_x R_W)^2)
 \\ &&- 2 \delta^{-\beta}  \partial_x (e^{\delta^{\alpha}W}(e^{\delta^{\beta}R_W}-1))+ \delta^{-\beta} Res_v,
\nonumber
\end{eqnarray*}
which we solve with $ R_W(\cdot,0) = R_V(\cdot,0) = 0 $.

\subsection{Estimates for the lowest order  derivatives}

\label{sec31}

We start  with the estimates for the lowest order derivatives. 
For proving our results this section is not  necessary, see Remark \ref{rem35}. However, we think
that it strongly helps to understand the calculations in the subsequent Section
\ref{secHigherEnergy}.
\begin{remark}\label{rem31}{\rm
The exponential terms are handled with a Taylor expansion. In detail
we set $z=\delta^\beta R_W$ and define $G$ by
$
\delta^\alpha G=e^{\delta^\alpha W}-1.
$
We use the identity
\begin{eqnarray*}
\lefteqn{\delta^{-\beta}\partial_x
\left(e^{\delta^\alpha W}(e^{\delta^\beta R_W}-1)\right)}\\
&=&\partial_xR_W+\delta^{1+\alpha}(\partial_XG)R_W
+\delta^\alpha G\partial_xR_W+\mathcal N, 
\end{eqnarray*}
where
\begin{eqnarray*}
\mathcal N&=&e^{\delta^\alpha W}\left(
(e^z-1)\partial_xR_W
+\delta^{1+\alpha}(\partial_XW)\delta^{-\beta}(e^z-1-z)\right).
\end{eqnarray*}
On any interval on which $\|\delta^\beta R_W\|_{L^\infty}\leq M$,
Taylor's formula for instance gives
\[
|e^z-1|\leq C_M|z|,\qquad
|e^z-1-z|\leq C_M|z|^2,\qquad
\left|e^z-1-z-\frac12z^2\right|\leq C_M|z|^3.
\]
Consequently, we find,
\[
|\mathcal N|\leq C_M\delta^\beta
\left(|R_W|\,|\partial_xR_W|+\delta^{1+\alpha}|R_W|^2\right).
\]}
\end{remark}
\medskip

The constants in the following estimates of this remainder may depend on $M$,
from Remark \ref{rem31},
and the  bounds on $ V $ and $ W $, but not on $\delta$.

The system for the error is of the form 
\begin{eqnarray*}
\partial_t R_W & = & - 2 \partial_x R_V + \ldots , \\
\partial_t R_V & = & -2 \partial_x R_W +  \frac12 \partial_x^3 R_W 
+ \ldots .
\end{eqnarray*}
To obtain cancellation of the constant-coefficient linear terms we consider the energy part
$$
\widetilde{\mathcal{E}}_0 = \int R_W^2  dx + \int R_V^2  dx + a \int ( \partial_x R_W)^2  dx.
$$ 
We find 
\begin{eqnarray*}
\frac12 \partial_t  \widetilde{\mathcal{E}}_0 & = & 
\int R_W \partial_t R_W  dx + \int R_V \partial_t R_V dx + a \int ( \partial_x R_W) ( \partial_x \partial_t R_W) dx \\ 
&= & 
\int R_W (- 2 \partial_x R_V )  dx + \int R_V (-2 \partial_x R_W +  \frac12 \partial_x^3 R_W)dx 
\\ && + a \int ( \partial_x R_W) ( \partial_x (- 2 \partial_x R_V))  dx + \ldots \\ 
 \\ 
&= & 
 \int R_V ( \frac12 \partial_x^3 R_W)dx  -2 a \int ( \partial_x^3 R_W)  R_V  dx + \ldots = 0
 + \ldots
\end{eqnarray*}
if we choose $a = \frac14 $. Before we estimate the nonlinear terms 
we introduce the auxiliary quantities
\[
E_m=\sum_{j=0}^m\left(\|\partial_x^jR_W\|_{L^2}^2
+\|\partial_x^jR_V\|_{L^2}^2\right),\qquad m\in\N_0.
\]
These quantities measure $H^m\times H^m$ regularity and are distinct from
$\mathcal E_0$, which controls $R_W$ in $H^1$ and $R_V$ in $L^2$.
For the residual estimates we write the approximation time scale as
$\delta^{-\tau}$, where
\[
\tau=\begin{cases}
3,&\text{for KdV and Airy},\\
1+\alpha,&\text{for Burgers},\\
1,&\text{for Whitham’s modulation}.
\end{cases}
\]
Young's inequality gives
\[
\delta^{-\beta}\|R\|_{L^2}\|Res\|_{L^2}
\leq\delta^\tau\|R\|_{L^2}^2
+\frac14\delta^{-(\tau+2\beta)}\|Res\|_{L^2}^2.
\]
As already said, we start with the lowest order derivatives. 
We find 
$$
\frac12 \partial_t \int R_W^2  dx = s_{a,1}+\ldots+ s_{a,5},
$$ 
where
\begin{eqnarray*}
 s_{a,1} & = & - 2 \int R_W \partial_x R_V dx , \\ 
 s_{a,2} & = &- 2 \delta^{\alpha} \int R_W V \partial_x R_W dx = \delta^{1+\alpha} \int R_W^2  \partial_X V dx, \\ 
 s_{a,3} & = &
 - 2 \delta^{1+  \alpha } \int R_W R_V \partial_X W dx, \\ 
 s_{a,4} & = &
 - 2 \delta^{\beta}  \int R_W R_V \partial_x R_W dx, \\ 
 s_{a,5} & = &
   \delta^{-\beta}\int R_W Res_w dx.
\end{eqnarray*}
We estimate
\begin{eqnarray*}
| s_{a,2} |+ | s_{a,3} | & \leq & C  \delta^{1+\alpha} E_0 , \\ 
| s_{a,4} | & \leq & 2 \delta^{\beta}  \|R_W \|_{C^0_b} \| R_V \|_{L^2} \| \partial_x R_W \|_{L^2} 
\leq C  \delta^{\beta} E_1^{3/2} , \\ 
| s_{a,5} | & \leq & \delta^{-\beta} \| R_W \|_{L^2} \| Res_w \|_{L^2}
\leq \delta^\tau \| R_W \|_{L^2}^2 + \delta^{-(\tau+2\beta)}  \| Res_w \|_{L^2}^2 \\ 
& \leq & C  \delta^\tau E_0 + \delta^{-(\tau+2\beta)}  \| Res_w \|_{L^2}^2.
\end{eqnarray*}
Moreover, we have 
$$
\frac18 \partial_t \int (\partial_x R_W)^2  dx = s_{a,6}+\ldots+ s_{a,10},
$$ 
where
\begin{eqnarray*}
 s_{a,6} & = & - \frac12 \int (\partial_x R_W) \partial_x^2 R_V dx , \\ 
 s_{a,7} & = &- \frac12 \delta^{\alpha} \int (\partial_x R_W) \partial_x(V \partial_x R_W) dx 
= - \frac14\delta^{1+\alpha} \int (\partial_x R_W)^2  \partial_X V dx
 , \\ 
 s_{a,8} & = &
 - \frac12 \delta^{1+  \alpha } \int (\partial_x R_W) \partial_x(R_V \partial_X W) dx, \\ 
 s_{a,9} & = &
 - \frac12 \delta^{\beta}  \int (\partial_x R_W) \partial_x(R_V \partial_x R_W) dx, \\ 
 s_{a,10} & = &
  \frac14 \delta^{-\beta}\int (\partial_x R_W) \partial_x Res_w dx.
\end{eqnarray*}
As above we estimate
\begin{eqnarray*}
| s_{a,7} |+ | s_{a,8} | & \leq & C  \delta^{1+\alpha} E_1 , \\ 
| s_{a,9} | & \leq &  C  \delta^{\beta} E_2^{3/2} , \\ 
| s_{a,10} | & \leq &  C  \delta^\tau E_1 + \delta^{-(\tau+2\beta)}  \| Res_w \|_{H^1}^2.
\end{eqnarray*}
Moreover, we find 
$$
 \frac12 \partial_t \int R_V^2 dx  =   s_{b,1}+\ldots+ s_{b,10},
$$
where
\begin{eqnarray*}
 s_{b,1} & = & \frac12 \int R_V \partial_x^3 R_W dx,\\
 s_{b,2} & = & - 2 \delta^{\alpha} \int  R_V  \partial_x (V R_V) dx  
 = - \delta^{1+ \alpha} \int  R_V^2  \partial_X V dx,
  \\ 
 s_{b,3} & = & 
 -  \delta^{\beta}  \int  R_V \partial_x (R_V^2)dx  = - \frac23 \delta^{\beta} 
 \int   \partial_x (R_V^3)dx = 0 ,
 \\ 
 s_{b,4} & = & 
\frac12 \delta^{1+ \alpha} \int  R_V \partial_x((\partial_X W)(\partial_x R_W))dx ,
 \\ 
 s_{b,5} & = & 
 \frac14\delta^{\beta} \int R_V \partial_x ((\partial_x R_W)^2)dx ,
 \\ 
 s_{b,6} & = & - 2 \int R_V 
 \partial_x R_W dx ,
 \\ 
 s_{b,7} & = & -2 \delta^{1+\alpha}\int R_V(\partial_X G) R_W
 dx ,
 \\ 
s_{b,8} & = & -2 \int R_V (\delta^{\alpha}G)  \partial_x R_W
dx ,
 \\ 
s_{b,9} & = & -2 \int R_V\mathcal N\,dx ,
 \\ 
 s_{b,10} & = & 
 \int R_V \delta^{-\beta} Res_v dx .
 \end{eqnarray*}
We first find that 
$$ 
s_{a,1} + s_{b,6} = 0  \qquad \textrm{and} \qquad s_{b,1} + s_{a,6} = 0 .
$$ 
We estimate
\begin{eqnarray*}
| s_{b,2} |+ | s_{b,4} | & \leq & C  \delta^{1+\alpha} E_2 , \\ 
| s_{b,5} | & \leq & 2 C \delta^{\beta}  \|R_V \|_{C_b^0} \| \partial_x R_W \|_{L^2} \| \partial_x^2 R_W \|_{L^2} 
\leq C  \delta^{\beta} E_2^{3/2} , \\ 
| s_{b,7} | & \leq & C  \delta^{1+\alpha} E_0 , \\ 
| s_{b,9} | & \leq & C_M\delta^\beta\|R_V\|_{L^2}\|R_W\|_{L^\infty}
\left(\|\partial_xR_W\|_{L^2}+\delta^{1+\alpha}\|R_W\|_{L^2}\right)
\\ & \leq & C_M\delta^\beta E_1^{3/2} , \\ 
| s_{b,10} | & \leq & \delta^{-\beta} \| R_V \|_{L^2} \| Res_v \|_{L^2}
\leq \delta^\tau \| R_V \|_{L^2}^2 + \delta^{-(\tau+2\beta)}  \| Res_v \|_{L^2}^2 \\ 
& \leq & C  \delta^\tau E_0 + \delta^{-(\tau+2\beta)}  \| Res_v \|_{L^2}^2.
\end{eqnarray*}
The coefficients $\delta^{1+\alpha}$ and $\delta^\beta$ have the required
size on the time scale $\delta^{-\tau}$ when $1+\alpha\geq\tau$ and
$\beta>\tau$, once the higher-derivative estimates are closed.
These inequalities hold for all approximations; in particular, for higher order Airy,
we have $\beta-\tau=\alpha-3/2>0$, and for higher order Burgers,
we have $\beta-\tau=1/2$.
Hence, it remains to get rid of  $ s_{b,8} $. In order to get rid of this term we extend the energy by a term similar to a term we already used in the derivation of the amplitude equations. 
We additionally consider the energy term
$b\delta^\alpha\int G R_W^2\,dx$, 
with $ G $ from Remark \ref{rem31}. Since $X=\delta(x-ct)$ and
$T=\delta^\tau t$, the chain rule gives
\[
\partial_tG(X,T)=-c\delta\partial_XG+\delta^\tau\partial_TG.
\]
We obtain
\begin{eqnarray*}
\lefteqn{\frac{b\delta^\alpha}{2}\partial_t\int G R_W^2\,dx} \\
&=&\frac{b\delta^\alpha}{2}\int
\left(-c\delta\partial_XG+\delta^\tau\partial_TG\right)R_W^2\,dx
+b\delta^\alpha\int G R_W\partial_tR_W\,dx\\
&=&\frac{b\delta^\alpha}{2}\int
\left(-c\delta\partial_XG+\delta^\tau\partial_TG\right)R_W^2\,dx
+\delta^\alpha\sum_{j=1}^5s_{c,j},
\end{eqnarray*}
where
\begin{eqnarray*}
 s_{c,1} & = & - 2 b  \int R_W G \partial_x R_V dx , \\ 
 s_{c,2} & = &- 2 b \delta^{\alpha} \int R_W G V \partial_x R_W dx = b \delta^{1+\alpha} \int R_W^2  \partial_X (GV) dx, \\ 
 s_{c,3} & = &
 - 2 b\delta^{1+  \alpha } \int R_W G R_V \partial_X W dx, \\ 
 s_{c,4} & = &
 - 2 b\delta^{\beta}  \int R_W G  R_V \partial_x R_W dx, \\ 
 s_{c,5} & = &
  b \delta^{-\beta}\int R_W G Res_w dx.
\end{eqnarray*}
We estimate as above
\begin{eqnarray*}
| s_{c,2} |+ | s_{c,3} | & \leq & C  \delta^{1+\alpha} E_0 , \\ 
| s_{c,4} | & \leq & C  \delta^{\beta} E_1^{3/2} , \\ 
| s_{c,5} | & \leq & 
 C  \delta^\tau E_0 + \delta^{-(\tau+2\beta)}  \| Res_w \|_{L^2}^2.
\end{eqnarray*}
The  derivative of $ G $ satisfies
\[
\left|\frac{b\delta^\alpha}{2}\int
\left(-c\delta\partial_XG+\delta^\tau\partial_TG\right)R_W^2\,dx\right|
\leq C\left(\delta^{1+\alpha}+\delta^{\alpha+\tau}\right)E_0,
\]
using uniform bounds on $\partial_XG$ and $\partial_TG$.
Finally, with $b=1$, we have
\begin{eqnarray*}
\lefteqn{\delta^{\alpha} s_{c,1} + s_{b,8}} \\  & = & - 2 b \delta^{\alpha} \int R_W G \partial_x R_V dx  -2 \int R_V (\delta^{\alpha}G)  \partial_x R_W
dx \\  & = & - 2 b \delta^{\alpha} \int R_W G \partial_x R_V dx  +2 \int R_V (\delta^{1+\alpha}\partial_X G)  R_W
dx \\ && +2 \int (\partial_x R_V) (\delta^{\alpha}G)   R_W
dx 
 \\  & = & 2 \delta^{1+\alpha} \int R_V (\partial_X G)  R_W
dx.
\end{eqnarray*}
Note that $\delta^\alpha s_{c,1}+s_{b,8}$ cancels $s_{b,7}$ exactly but we will not make use of this property.
With $a=1/4$ and $b=1$, the resulting energy is
\begin{eqnarray} \label{preenergy}
\mathcal E& =& \int\left(R_W^2+\frac14(\partial_xR_W)^2+R_V^2
+\delta^\alpha G R_W^2\right)\,dx \\
& = & \int\left(e^{\delta^\alpha W}R_W^2+R_V^2
+\frac14(\partial_xR_W)^2\right)\,dx. \nonumber
\end{eqnarray}
\begin{remark}{\rm  \label{remequi}
Uniform boundedness of $ W$ implies that $\mathcal E$ is
uniformly equivalent to $\|R_W\|_{H^1}^2+\|R_V\|_{L^2}^2$ for $ \delta > 0 $ sufficiently small if $ \alpha > 0 $. 
For $ \alpha = 0 $, if we use the first line of \eqref{preenergy},
it seems that  an additional smallness condition 
on $ C_1 $ defined in $ \eqref{condition1} $ 
in Theorem \ref{th11whit} is needed.
However, the representation of \eqref{preenergy} in the second line 
shows that this is not necessary.
Moreover, as already said, Section \ref{sec31} is not necessary
for the proof  of the above approximation results. See Remark \ref{rem35}.
For  the modified energy, used in the subsequent
Section \ref{secHigherEnergy}, such a smallness condition also does not occur.
See \eqref{higherWeightedEnergy} and \eqref{higherEnergyEquivalence}.
}
\end{remark}
\begin{remark}{\rm  \label{rem32}
To estimate the residual contributions, we require the  bound
\begin{equation} \label{rem32b}
\sup_{t\in[0,T_0/\delta^\tau]}
\|(Res_v,Res_w)(\cdot,t)\|_{L^2\times H^1}
\leq C\delta^{\beta+\tau}.
\end{equation}
The required residual powers are $13/2$ for KdV, $\alpha+9/2$ for
higher order Airy, $2\alpha+5/2$ for higher order Burgers, and $5/2$
for Whitham’s modulation.
These powers  are established in Lemma \ref{lemreskdv}-Lemma \ref{lemreswme}.
}
\end{remark}
The bound \eqref{rem32b} implies that
\[
\delta^{-(\tau+2\beta)}
\|(Res_v,Res_w)\|_{L^2\times H^1}^2\leq C\delta^\tau.
\]
and so combining the preceding estimates and absorbing fixed factors into the
constants gives
\begin{equation} \label{errorlow}
\partial_t\mathcal E
\leq C(\delta^{1+\alpha}+\delta^\tau)E_2
+C\delta^\beta E_2^{3/2}+C\delta^\tau.
\end{equation}
All terms on the right hand side have at least a $ \delta^{\tau} $ in front.
Therefore, in principle this inequality allows us to estimate the error on the long   
$ \delta^{-\tau} $ time scale. However, the time derivative of the quantity $  \mathcal E $
which contains terms with at most one derivative
is estimated by quantities with up to two derivatives.
Hence,  \eqref{errorlow} is not sufficient to bound $  \mathcal E $ on the long   
$ \delta^{-\tau} $ time scale by Gronwall's inequality.

\subsection{Estimates for the highest derivatives}

\label{secHigherEnergy}

The estimate \eqref{errorlow} will now be closed with the associated estimates for  the highest derivatives, i.e., in contrast to Section \ref{sec31} the right hand side 
of an inequality similar to \eqref{errorlow} 
should not contain more spatial derivatives than the left hand side.

\paragraph{Notation.}
For the subsequent calculations we use the following  notation.
For the approximation we write 
$$ 
q_W=\delta^\alpha W(X,T),\qquad q_V=\delta^\alpha V(X,T).
$$ 
We write the exact  
solution as 
$$ 
w=q_W+\delta^{\beta} R_W,\qquad v=q_V+\delta^{\beta} R_V .
$$
For closing the estimates we also need the quantity 
$$ 
R_Z=\partial_x R_W.
$$ 
The rescaled residual terms are abbreviated by 
$$f_w=\delta^{-\beta}Res_w,\qquad f_v=\delta^{-\beta}Res_v.$$
We introduce the abbreviation 
$$ 
G_1(R_W)=\frac{e^{\delta^{\beta} R_W}-1}{\delta^{\beta}}
=R_W \int_0^1e^{\sigma\delta^{\beta} R_W}\,d\sigma
$$
and the operator 
$$
L_w=\partial_x^2+(\partial_x w)\partial_x
=e^{- w}\partial_x(e^w \partial_x).
$$
Finally, the derivatives are abbreviated as  
$$ R_{W,k}=\partial_x^k R_W,  \qquad R_{V,k}=\partial_x^k R_{V},  \qquad R_{Z,k}=\partial_x^k R_{Z}.$$

\paragraph{Rewriting the error equations.}
With these notations and abbreviations the equation for the error function 
$ R_W $ is  written as 
\begin{equation}
\partial_t R_W=-2\partial_x R_V-2v\partial_x R_W-2(\partial_x q_W)R_V+f_w,
\label{higherErrorR}
\end{equation}
Differentiating \eqref{higherErrorR}
with respect to $x$ gives
\begin{align} \nonumber
\partial_t R_Z={}&-2\partial_x^2R_V-2v\partial_x R_Z
-2(\partial_x q_V+\delta^{\beta}\partial_x R_V)R_Z\\
&-2(\partial_x^2q_W)R_V-2(\partial_x q_W)\partial_x R_V
+\partial_xf_w \nonumber
\\ = & -2L_w R_V-2v\partial_x R_Z
-2(\partial_x q_V)R_Z-2(\partial_x^2q_W)R_V+\partial_xf_w,
\label{higherErrorP}
\end{align}
where we used 
$$ 
-2L_w R_V = -2\partial_x^2R_V-2(\partial_x q_W+\delta^{\beta} R_Z)\partial_x R_V.
$$ 
Finally the equation for the error function 
$ R_V $ is  written as
\begin{align}
 \label{higherErrorZ}
\partial_t R_V&=\frac12L_w R_Z-2v\partial_x R_V-2e^w R_Z
+\frac12(\partial_x^2q_W)R_Z-2(\partial_x q_V)R_V
\\ & \qquad -2(\partial_x q_W)e^{q_W} G_1(R_W)+f_v.
\nonumber
\end{align}
where we used 
\[
\frac12\partial_x^3R_W
+\frac12\partial_x\bigl((\partial_xq_W)R_Z\bigr)
+\frac{\delta^{\beta}}{4}\partial_x(R_Z^2)
=\frac12L_wR_Z+\frac12(\partial_x^2q_W)R_Z,
\]
\[
-2 \delta^{-\beta}\partial_x(e^w-e^{q_W})
=-2 e^w R_Z-2(\partial_x q_W)e^{q_W}G_1(R_W),
\]
and 
\[
-\delta^{-\beta}\partial_x(v^2-{q_V}^2)
=-2v\partial_x R_V-2(\partial_x q_V)R_V.
\]

\paragraph{Introducing the energy.}
For estimating the $ k $-th derivatives $ R_{W,k} $, $ R_{V,k} $, and $ R_{Z,k} $ of the error functions $ R_{W} $, $ R_{V} $, and $ R_{Z} $,
for $0\leq k\leq m$, we define the energies 
\begin{equation}\label{higherWeightedEnergy}
\mathcal E_k=\int\left(e^{2 w}R_{W,k}^2
+e^w R_{V,k}^2+\frac14e^w R_{Z,k}^2\right)\,dx,
\qquad \mathcal E^{(m)}=\sum_{k=0}^m\mathcal E_k.
\end{equation}
The full-solution weights in $\mathcal E_0$ are chosen to cancel
the highest-order couplings. In particular,
\[
\int e^w fL_wg\,dx=-\int e^w(\partial_xf)(\partial_xg)\,dx.
\]
The common weight $e^w$ 
cancels the principal second-order coupling between $R_V$ and
$R_Z=\partial_xR_W$. The cancelation by the weight $e^{2w}$ on $R_W$ 
allows the higher-derivative
estimates to close without derivative loss.

If $ w $ is uniformly bounded, then
\begin{equation}\label{higherEnergyEquivalence}
c\bigl(\|R_{W}\|_{H^{m+1}}^2+\|R_{V}\|_{H^m}^2\bigr)
\leq\mathcal E^{(m)}
\leq C\bigl(\|R_{W}\|_{H^{m+1}}^2+\|R_{V}\|_{H^m}^2\bigr).
\end{equation}

\paragraph{Differentiating the energy terms:}

{\bf a)} Differentiating \eqref{higherErrorR} $k$ times gives
\begin{align*}
\partial_t R_{W,k}={}&-2\partial_x R_{V,k}-2v\partial_x R_{W,k}
-2\sum_{j=1}^k\binom{k}{j}(\partial_x^jv)R_{W,k-j+1}\\
&-2\sum_{j=0}^k\binom{k}{j}(\partial_x^{j+1}q_W)R_{V,k-j}
+\partial_x^kf_w.
\end{align*}
The product rule yields
\[
\frac{d}{dt}\int e^{2w}R_{W,k}^2\,dx
=2\int e^{2w}(\partial_tw)R_{W,k}^2\,dx
+2\int e^{2w}R_{W,k}\partial_t R_{W,k}\,dx.
\]
To rewrite the first term  on the right hand side we use the original equation  
\begin{equation}\label{higherWeightEvolution}
\partial_t w=-2\partial_x v-2v\partial_x w
\end{equation}
such that
\begin{align*}
2\int e^{2w}(\partial_t w)R_{W,k}^2\,dx
={}&-4\int e^{2w}(\partial_x v)R_{W,k}^2\,dx
-4\int e^{2w}v(\partial_x w)R_{W,k}^2\,dx.
\end{align*}
Since $\partial_x(e^{2w})=2e^{2w}\partial_x w$, integration
by parts in the second integral yields
\begin{align*}
-4\int e^{2w}v(\partial_x w)R_{W,k}^2\,dx
&=-2\int \partial_x(e^{2w})vR_{W,k}^2\,dx\\
&=2\int e^{2w}\partial_x\bigl(vR_{W,k}^2\bigr)\,dx\\
&=2\int e^{2w}(\partial_x v)R_{W,k}^2\,dx
+4\int e^{2w}vR_{W,k}\partial_x R_{W,k}\,dx.
\end{align*}
Thus
\begin{eqnarray*}
2\int e^{2w}(\partial_t w)R_{W,k}^2\,dx
& = & 4\int e^{2w}vR_{W,k}\partial_x R_{W,k}\,dx
-2\int e^{2w}(\partial_x v)R_{W,k}^2\,dx \\ & = & s_{W,1} + s_{W,2}.
\end{eqnarray*}
Therefore,  the full derivative is
\[
\frac{d}{dt}\int e^{2w}R_{W,k}^2\,dx
=s_{W,1}+\cdots+s_{W,7},
\]
where \begin{align*}
s_{W,1}&= 4\int e^{2w}vR_{W,k}\partial_x R_{W,k}\,dx ,\\
s_{W,2}&=-2\int e^{2w}(\partial_x v)R_{W,k}^2\,dx ,\\
s_{W,3}&=-4\int e^{2w}R_{W,k} \partial_x R_{V,k}\,dx,\\
s_{W,4}&=-4\int e^{2w}R_{W,k} v\partial_x R_{W,k}\,dx,\\
s_{W,5}&=-4\int e^{2w}R_{W,k} \sum_{j=1}^k\binom{k}{j}(\partial_x^jv)R_{W,k-j+1}\,dx,\\
s_{W,6}&=-4\int e^{2w}R_{W,k} \sum_{j=0}^k\binom{k}{j}(\partial_x^{j+1}q_W)R_{V,k-j}\,dx,\\
s_{W,7}&=2\int e^{2w}R_{W,k}\partial_x^kf_w dx.
\end{align*}
For $k=0$, the term $ s_{W,5} $ vanishes. 
We see that 
$$ 
s_{W,1} + s_{W,4} = 0 .
$$ 
We rewrite $ s_{W,3} $, again by 
integration by parts, into
\begin{align*}
-4\int e^{2w}R_{W,k}\partial_x R_{V,k}\,dx
&=4\int \partial_x\bigl(e^{2w}R_{W,k}\bigr)R_{V,k}\,dx\\
&=4\int \bigl(e^{2w}\partial_x R_{W,k}
+2e^{2w}(\partial_x w)R_{W,k}\bigr)R_{V,k}\,dx\\
&=4\int e^{2w}R_{Z,k}R_{V,k}\,dx
+8\int e^{2w}(\partial_x w)R_{W,k}R_{V,k}\,dx
\\ & = s_{W,8} + s_{W,9}
\end{align*}
where we used $\partial_x R_{W,k}=R_{Z,k}$.
Therefore, we finally have 
\begin{equation}
\frac{d}{dt}\int e^{2w}R_{W,k}^2\,dx
=s_{W,2}+ s_{W,5}+ s_{W,6}+s_{W,7}+ s_{W,8}+s_{W,9}.
\end{equation}
\medskip

{\bf b)} Differentiating \eqref{higherErrorZ} $k$ times gives
\begin{align*}
\partial_t R_{V,k}={}&\frac12L_w R_{Z,k}
+\frac12\sum_{j=1}^k\binom{k}{j}
(\partial_x^{j+1}w)R_{Z,k-j+1}
-2v\partial_x R_{V,k}\\
&-2\sum_{j=1}^k\binom{k}{j}(\partial_x^jv)R_{V,k-j+1}
-2e^w R_{Z,k}
-2\sum_{j=1}^k\binom{k}{j}(\partial_x^je^w)R_{Z,k-j}\\
&+\frac12(\partial_x^2q_W)R_{Z,k}
+\frac12\sum_{j=1}^k\binom{k}{j}(\partial_x^{j+2}q_W)R_{Z,k-j}
-2(\partial_x q_V)R_{V,k}\\
&-2\sum_{j=1}^k\binom{k}{j}(\partial_x^{j+1}q_V)R_{V,k-j}
-2\partial_x^k\bigl((\partial_x q_W)e^{q_W}G_1(R_W)\bigr)
+\partial_x^kf_v.
\end{align*}
The
product rule yields
\[
\frac{d}{dt}\int e^w R_{V,k}^2\,dx
=\int e^w(\partial_tw)R_{V,k}^2\,dx
+2\int e^w R_{V,k}\partial_t R_{V,k}\,dx.
\]
Again by \eqref{higherWeightEvolution} we obtain  
for the first term on the right hand side that 
\begin{eqnarray*}
\int e^w(\partial_tw)R_{V,k}^2\,dx & = & -2 \int e^w(\partial_x v)R_{V,k}^2\,dx -2 \int e^w(v \partial_xw)R_{V,k}^2\,dx \\ 
& = & -2 \int \partial_x (e^w v)R_{V,k}^2\,dx \\
& = & 4\int e^w v R_{V,k}\partial_x R_{V,k}\,dx = s_{V,1}.
\end{eqnarray*}
Therefore, we have 
\[
\frac{d}{dt}\int e^w R_{V,k}^2\,dx
=s_{V,1}+\cdots+s_{V,13},
\]
where 
\begin{align*}
s_{V,1} & = 4\int e^w v R_{V,k}\partial_x R_{V,k}\,dx , \\
s_{V,2} & = \int e^w R_{V,k} L_w R_{Z,k},dx 
 =\int R_{V,k}\partial_x(e^w\partial_x R_{Z,k})\,dx \\
&=-\int e^w(\partial_x R_{V,k})(\partial_x R_{Z,k})\,dx
,\\
s_{V,3} & = \int e^w R_{V,k}\sum_{j=1}^k\binom{k}{j} (\partial_x^{j+1}w)R_{Z,k-j+1}dx,\\
s_{V,4} & = -4 \int e^w R_{V,k}v\partial_x R_{V,k},dx,\\
s_{V,5} & = -4\int e^w R_{V,k}\sum_{j=1}^k\binom{k}{j}(\partial_x^jv)R_{V,k-j+1}dx,\\
s_{V,6} & = -4\int e^w R_{V,k}e^w R_{Z,k},dx,\\
s_{V,7} & = -4\int e^w R_{V,k}\sum_{j=1}^k\binom{k}{j}(\partial_x^je^w)R_{Z,k-j}dx,\\
s_{V,8} & = \int e^w R_{V,k}(\partial_x^2q_W)R_{Z,k},dx,\\
s_{V,9} & = \int e^w R_{V,k}\sum_{j=1}^k\binom{k}{j}(\partial_x^{j+2}q_W)R_{Z,k-j}dx,\\
s_{V,10} & = -4\int e^w R_{V,k}(\partial_x q_V)R_{V,k},dx,\\
s_{V,11} & = -4\int e^w R_{V,k}\sum_{j=1}^k\binom{k}{j}(\partial_x^{j+1}q_V)R_{V,k-j}dx,\\
s_{V,12} & = -4\int e^w R_{V,k}\partial_x^k\bigl((\partial_x q_W)e^{q_W}G_1(R_W)\bigr)dx,\\
s_{V,13} & = 2\int e^w R_{V,k} \partial_x^kf_vdx.
\end{align*}
We see that 
$$ 
s_{V,1} + s_{V,4} = 0 .
$$ 
For $k=0$, the terms $ s_{V,3} $, $ s_{V,5} $, $ s_{V,7} $, $ s_{V,9} $, and $ s_{V,11} $ vanish. Again all integrations by parts above have vanishing boundary terms.
\medskip

{\bf c)} Differentiating \eqref{higherErrorP} $k$ times gives
\begin{align*}
\partial_t R_{Z,k}={}&-2L_w R_{V,k}
-2\sum_{j=1}^k\binom{k}{j}
(\partial_x^{j+1}w)R_{V,k-j+1}
-2v\partial_x R_{Z,k}\\
&-2\sum_{j=1}^k\binom{k}{j}(\partial_x^jv)R_{Z,k-j+1}
-2(\partial_x q_V)R_{Z,k}
-2\sum_{j=1}^k\binom{k}{j}(\partial_x^{j+1}q_V)R_{Z,k-j}\\
&-2(\partial_x^2q_W)R_{V,k}
-2\sum_{j=1}^k\binom{k}{j}(\partial_x^{j+2}q_W)R_{V,k-j}
+\partial_x^{k+1}f_w.
\end{align*}
The product rule yields
\[
\frac{d}{dt}\int\frac14e^w R_{Z,k}^2\,dx
=\frac14\int e^w(\partial_tw)R_{Z,k}^2\,dx
+\frac12\int e^w R_{Z,k}\partial_t R_{Z,k}\,dx.
\]
Again by \eqref{higherWeightEvolution} we obtain  
for the first term on the right hand side that 
\begin{eqnarray*}
\frac14\int e^w(\partial_tw)R_{Z,k}^2\,dx & = & -\frac12 \int e^w(\partial_x v)R_{Z,k}^2\,dx -\frac12 \int e^w(v \partial_xw)R_{Z,k}^2\,dx \\ 
& = & -\frac12 \int \partial_x (e^w v)R_{Z,k}^2\,dx \\
& = & \int e^w v R_{Z,k}\partial_x R_{Z,k}\,dx = s_{Z,1}.
\end{eqnarray*}
Therefore, we have 
\[
\frac{d}{dt}\int\frac14e^w R_{Z,k}^2\,dx
=s_{Z,1}+\cdots+s_{Z,10},
\]
where 
\begin{align*}
s_{Z,1} & =  \int e^w v R_{Z,k}\partial_x R_{Z,k}\,dx, \\
s_{Z,2} & =  -\int e^w  R_{Z,k}L_w R_{V,k}\,dx 
=-\int R_{Z,k}\partial_x(e^w\partial_x R_{V,k})\,dx \\
& =\int e^w(\partial_x R_{Z,k})(\partial_x R_{V,k})\,dx,\\
s_{Z,3} & =  - \int e^w  R_{Z,k}\sum_{j=1}^k\binom{k}{j} (\partial_x^{j+1}w)R_{V,k-j+1}\,dx, \\
s_{Z,4} & =  - \int e^w  R_{Z,k}v\partial_x R_{Z,k}\,dx, \\
s_{Z,5} & =  - \int e^w  R_{Z,k}\sum_{j=1}^k\binom{k}{j}(\partial_x^jv)R_{Z,k-j+1}\,dx, \\
s_{Z,6} & =  -\int e^w  R_{Z,k}(\partial_x q_V)R_{Z,k}\,dx, \\
s_{Z,7} & =  -\int e^w  R_{Z,k}\sum_{j=1}^k\binom{k}{j}(\partial_x^{j+1}q_V)R_{Z,k-j}\,dx, \\
s_{Z,8} & =  -\int e^w  R_{Z,k}(\partial_x^2q_W)R_{V,k}\,dx, \\
s_{Z,9} & =  -\int e^w  R_{Z,k}\sum_{j=1}^k\binom{k}{j}(\partial_x^{j+2}q_W)R_{V,k-j}\,dx, \\
s_{Z,10} & =  \frac12 \int e^w  R_{Z,k}\partial_x^{k+1}f_w\,dx.
\end{align*}
For $k=0$, the terms $ s_{Z,3} $, $ s_{Z,5} $, $ s_{Z,7} $,  and $ s_{Z,9} $ vanish. Again all integrations by parts above have vanishing boundary terms.

\paragraph{Cancelations.}
Additional to the cancelations 
$$ 
s_{W,1} + s_{W,4} = 0 , \qquad s_{V,1} + s_{V,4} = 0 ,
$$ 
which we already found 
there is  the cancelation 
$$
s_{V,2}+ s_{Z,2} = 0 .
$$ 
Moreover, 
we have
\[
s_{Z,1}+s_{Z,4}=0,\qquad
s_{W,8}+s_{V,6}=0,\qquad
s_{V,8}+s_{Z,8}=0.
\]
After the preceding cancellations, the derivative of $\mathcal E_k$
contains exactly the following terms:
\begin{align*}
\frac{d}{dt}\mathcal E_k={}&s_{W,2}+s_{W,5}+s_{W,6}+s_{W,7}+s_{W,9}\\
&+s_{V,3}+s_{V,5}+s_{V,7}+s_{V,9}
+s_{V,10}+s_{V,11}+s_{V,12}+s_{V,13}\\
&+s_{Z,3}+s_{Z,5}+s_{Z,6}+s_{Z,7}+s_{Z,9}+s_{Z,10}.
\end{align*}
In the following we obtain bounds for the right-hand side of this inequality 
in terms of the full energy $ \mathcal E^{(m)} $.

\paragraph{Bounds for the approximation and the solution.}
Let $ m \in \N $ with $m\geq2$ and assume 
for $1\leq j\leq m+2$ and $0<\delta\leq1$ the uniform bounds
\begin{equation}\label{higherProfileBounds}
\|q_W\|_{L^\infty}+\|q_V\|_{L^\infty}\leq C_A,
\qquad
\|\partial_x^jq_W\|_{L^\infty}
+\|\partial_x^jq_V\|_{L^\infty}
\leq C_A\delta^{1+\alpha}.
\end{equation}
The last  bounds hold true due to
$\partial_x=\delta\partial_X$.
For closing the subsequent arguments  
we only assumed $ \leq C_A $ in the 
first bound in \eqref{higherProfileBounds} although $ \leq C_A \delta^{\alpha} $ seems 
more plausible.
Then as long as $ \sqrt{\mathcal E^{(m)}}\leq M$ for a $ M = \mathcal{O}(1) $, 
suitably chosen below, 
we have the estimates  
\begin{equation} \label{bound1}
\|\partial_x w\|_{L^\infty}+\|\partial_x v\|_{L^\infty}
+\|\partial_t w\|_{L^\infty}
\leq C_M\left(\delta^{1+\alpha}
+\delta^{\beta}\sqrt{\mathcal E^{(m)}}\right)
\end{equation}
due to $w=q_W+\delta^{\beta} R_W$ and $v=q_V+\delta^{\beta} R_V$,
where the time derivative is estimated using \eqref{higherWeightEvolution}.
Here and below, constants may also depend on $m$ and the approximation bounds
$ C_A $.

\paragraph{Bounds for $s_{W,2}$ and $s_{W,9}$.}
The preceding  bound \eqref{bound1} and the Cauchy-Schwarz inequality imply
\[
|s_{W,2}|+|s_{W,9}|
\leq C_M\left(\delta^{1+\alpha}
+\delta^{\beta}\sqrt{\mathcal E^{(m)}}\right)\mathcal E^{(m)}.
\]

\paragraph{Bounds for $ s_{W,5} $, $ s_{V,3}$, $s_{V,5}$, $s_{V,7}$,
$ s_{Z,3}$, and $s_{Z,5}$.}
In all these terms at least one derivative falls on $ v $ or $ w $ such that these terms 
are at least of order $ \mathcal{O}(\delta^{1+\alpha}) $. We have to take care 
that at most $ k $ derivatives fall on $ v $,  $ R_V $ and $ R_Z $ and 
at most $ k+1 $ derivatives fall on $ w $ and $ R_W $. 
This is satisfied for all these terms.

\[
|s_{W,5}|+|s_{V,3}|+|s_{V,5}|+|s_{V,7}|
+|s_{Z,3}|+|s_{Z,5}|
\leq C_M\left(\delta^{1+\alpha}
+\delta^{\beta}\sqrt{\mathcal E^{(m)}}\right)\mathcal E^{(m)}.
\]
All these  terms vanish for $k=0$.

\paragraph{Bounds for
$s_{W,6}$, $s_{V,9}$, $s_{V,10}$, $s_{V,11}$,
$s_{Z,6}$, $s_{Z,7}$, and $s_{Z,9}$.}
The  bounds on the approximation  give
\begin{align*}
&|s_{W,6}|+|s_{V,9}|+|s_{V,10}|+|s_{V,11}|\\
&\qquad+|s_{Z,6}|+|s_{Z,7}|+|s_{Z,9}|
\leq C_M\delta^{1+\alpha}\mathcal E^{(m)}.
\end{align*}

\paragraph{Bounds for $s_{V,12}$.}
With
\[
G_1(R_W)=R_W\int_0^1e^{\sigma\delta^{\beta} R_W}\,d\sigma
\]
we obtain
\[
\left\|\partial_x^k\bigl((\partial_x q_W)e^{q_W}
G_1(R_W)\bigr)\right\|_{L^2}
\leq C_M\delta^{1+\alpha}\sqrt{\mathcal E^{(m)}},
\qquad
|s_{V,12}|\leq C_M\delta^{1+\alpha}\mathcal E^{(m)}.
\]
\begin{remark}{\rm
Summing over $0\leq k\leq m$, all remaining non-forcing terms are
bounded in absolute value by
\begin{equation}\label{higherNonlinearBound}
C_M\delta^{1+\alpha}\mathcal E^{(m)}
+C_M\delta^\beta\bigl(\mathcal E^{(m)}\bigr)^{3/2}.
\end{equation}}
\end{remark}

\paragraph{Bounds for  $s_{W,7}$, $s_{V,13}$, and $s_{Z,10}$.}
We define
\[
C_{res,m}(t)=\|(Res_v,Res_w)(\cdot,t)\|_{H^m\times H^{m+1}}.
\]
Since $f_w=\delta^{-\beta}Res_w$ and $f_v=\delta^{-\beta}Res_v$, the
Cauchy-Schwarz inequality and Young's inequality immediately give
\begin{align*}
\sum_{k=0}^m\bigl(|s_{W,7}|+|s_{V,13}|+|s_{Z,10}|\bigr)
&\leq C_M\delta^{-\beta}\sqrt{\mathcal E^{(m)}}\,C_{res,m}\\
&\leq C_M\delta^\tau\mathcal E^{(m)}
+C_M\delta^{-(\tau+2\beta)}C_{res,m}^2.
\end{align*}
The $H^{m+1}$-norm of $Res_w$ is required by $s_{Z,10}$,
which contains $\partial_x^{k+1}f_w$.

\paragraph{Application of Gronwall's inequality}
Combining the preceding estimates gives
\begin{equation}\label{higherClosedEnergy}
\frac{d}{dt}\mathcal E^{(m)}
\leq C_1(\delta^{1+\alpha}+\delta^\tau)\mathcal E^{(m)}
+C_M\delta^\beta\bigl(\mathcal E^{(m)}\bigr)^{3/2}
+C_2\delta^{-(\tau+2\beta)}C_{res,m}^2,
\end{equation}
where $ C_2 $ can be chosen  independent of 
$ M $ and $ 0 \leq \delta  \ll 1 $.
The constant $ C_1 $ can be chosen  independent of 
$ M $ and $ 0 < \delta  \ll 1 $ if $ \delta^{\beta-\tau} M > 0 $ is sufficiently small.
This inequality contains no higher error norm on its right-hand side.
Since we have 
\[
\sup_{t\in[0,T_0/\delta^\tau]}C_{res,m}(t)
\leq C_A\delta^{\beta+\tau},
\]
for all approximations and $1+\alpha\geq\tau$
and $\beta>\tau$, we find that  \eqref{higherClosedEnergy} reduces to
\[
\frac{d}{dt}\mathcal E^{(m)}
\leq \delta^\tau(C_1 \mathcal E^{(m)}+C_2 C_A^2)
+C_M\delta^\beta\bigl(\mathcal E^{(m)}\bigr)^{3/2}.
\]
For 
$$ 
C_M \delta^{\beta-\tau} \bigl(\mathcal E^{(m)}\bigr)^{3/2}  \leq 1
$$ 
we have 
$$ 
\frac{d}{dt}\mathcal E^{(m)}
\leq \delta^\tau(C_1 \mathcal E^{(m)}+C_2 C_A^2+ 1).
$$ 
Gronwall's inequality yields 
$$ 
\mathcal E^{(m)} \leq (C_2 C_A^2+ 1) T_0 e^{C_1 T_0} =: M^2.
$$ 
Choosing $ \delta_0 > 0 $ so small that 
$$ 
C_M \delta_0^{\beta-\tau} M^{3}  \leq 1,
$$ 
we are done. \qed
\begin{remark}\label{weights}{\rm
The energy $\mathcal E_0$, defined in \eqref{higherWeightedEnergy},
is a nonlinear weighted refinement of
the energy $\mathcal E$,  defined in \eqref{preenergy}. Writing
\[
q_W=\delta^\alpha W,\qquad
w=q_W+\delta^\beta R_W,\qquad
B=R_V^2+\frac14(\partial_xR_W)^2,
\]
we have
\[
\mathcal E=\int(e^{q_W}R_W^2+B)\,dx,
\qquad
\mathcal E_0=\int(e^{2w}R_W^2+e^wB)\,dx.
\]
Consequently,
\[
\mathcal E_0-\mathcal E
=\int\left[(e^{2w}-e^{q_W})R_W^2+(e^w-1)B\right]\,dx.
\]
For uniformly bounded $q_W$ and $w$, this implies
\[
|\mathcal E_0-\mathcal E|
\leq C\left(\delta^\alpha\|W\|_{L^\infty}
+\delta^\beta\|R_W\|_{L^\infty}\right)\mathcal E.
\]
Thus, for $\alpha,\beta>0$ and uniformly bounded approximations and errors,
both energies have the same leading quadratic part,
\[
\int\left(R_W^2+R_V^2+\frac14(\partial_xR_W)^2\right)\,dx.
\]
However, it cannot be expected that the terms of order $ \mathcal{O}(\delta^{\alpha}) $ coincide since there are highest order derivative terms 
of this order 
which are handled differently in Section \ref{sec31} and Section \ref{secHigherEnergy}.
}\end{remark}
\begin{remark} \label{rem35}{\rm
For sufficiently regular solutions, the estimates of Section \ref{secHigherEnergy}
also yield the bounds for the lower-order derivatives which we found 
in Section \ref{sec31}.
Indeed, for $k=0$, as mentioned above, many terms vanish, and so we have
\begin{align*}
\frac{d}{dt}\mathcal E_0
\leq{}&
C_M\left(\delta^{1+\alpha}+\delta^\tau
+\delta^\beta\sqrt{\mathcal E^{(m)}}\right)\mathcal E_0\\
&+C_M\delta^{-(\tau+2\beta)}
\|(Res_v,Res_w)\|_{L^2\times H^1}^{\,2}.
\end{align*}
Here $m\geq2$ provides the required bounds on the first spatial
derivatives of the full solution. The uniform higher-order
estimate therefore implies $ 
\mathcal E_0(t)\leq C $.
Since  the energies satisfy
\[
\mathcal E\asymp\mathcal E_0
\asymp\|R_W\|_{H^1}^2+\|R_V\|_{L^2}^2,
\]
we see that Section \ref{sec31} is not needed to obtain this bound when the
higher regularity required in Section \ref{secHigherEnergy} is available.
}\end{remark}

\section{Discussion}

\label{secdisc}

We close this paper firstly with a remark on the general long wave limit which leads to decoupled counter-propagating waves  for 
$ \alpha > 0 $.
The second remark is about the reconstruction of $ \phi $ from $ w $ and $ v $.
\begin{remark}{\rm
The one-directional approximations considered above are part of a more
general long-wave description. For sufficiently regular, spatially
localized small-amplitude initial data, the variables $ (w,v) $ of both equations are
usually excited. The solution splits into a wave packet moving to the
right and a wave packet moving to the left. Indeed, the linearization
of \eqref{made1}-\eqref{made2} at $(w,v)=(0,0)$ has long-wave
characteristic speeds $2$ and $-2$. Therefore, for $ \alpha>0 $ the leading order approximation 
 has the form
\[
\begin{pmatrix}w\\v\end{pmatrix}(x,t)
\approx\delta^\alpha\left[
A_+(\delta(x-2t),T)\begin{pmatrix}1\\1\end{pmatrix}
+A_-(\delta(x+2t),T) \begin{pmatrix}1\\-1\end{pmatrix}
\right],
\]
with the slow time $T$ chosen according to the scaling regime.
In the KdV regime, each packet evolves to leading order according to
its own KdV equation, with the corresponding propagation direction.
The coupling terms between the two packets contribute only to
higher-order corrections to this leading order approximation.
Such a decoupling has been justified for instance for a
Boussinesq equation in \cite{Schn98LW} and for the water-wave
problem  in \cite{SW00,SW02}.
A corresponding two-directional justification for the present system requires estimates for these interaction terms. As far as we can see, this is rather straightforward.
}\end{remark}

\begin{remark}{\rm
For a sufficiently smooth solution $(v,w)$ of
\eqref{made1}-\eqref{made2} on $\mathbb R$, set $\rho=e^w>0$.
Fix $x_0\in\mathbb R$ and define
\[
\phi(x,t)=\phi_0(t)+\int_{x_0}^x v(y,t)\,dy.
\]
Choose the time-dependent integration
constant by
\[
\phi_0'(t)=\left[\frac12\partial_x^2w
+\frac14(\partial_xw)^2-v^2-2(e^w-1)\right]_{x=x_0}.
\]
See \cite{BKZ21} for more details.
}\end{remark}

\appendix

\section{Higher order Airy approximation}
\label{appHigherAiry}

This section contains a proof of Lemma
\ref{lemresairy} for the case $\alpha=3$. For this value we 
exemplarily construct a higher order Airy
approximation
 with residual $ \mathcal{O}(\delta^{15/2})
= \mathcal{O}(\delta^{\alpha+9/2})$ in $H^m\times H^{m+1}$.
We make the ansatz
\[
w_{\rm app}(x,t)=\delta^3W(X,T),\qquad
v_{\rm app}(x,t)=\delta^3V(X,T),
\]
with $ X=\delta(x-2t) $, $ T=\delta^3t $ and 
\[
W=W_0+\delta W_1,\qquad
V=W_0+\delta W_1+\delta^2W_2+\delta^3W_3,
\]
where the  $W_j$ are independent of $0 < \delta  \ll 1$. Since
$\partial_t=-2\delta\partial_X+\delta^3\partial_T$ and
$\partial_x=\delta\partial_X$, substitution into
\eqref{made1}-\eqref{made2} gives the residuals
\begin{align*}
Res_w={}&2\delta^4\partial_X(W-V)
-\delta^6\partial_TW-2\delta^7V\partial_XW,\\
Res_v={}&2\delta^4\partial_X(V-W)
-\delta^6\partial_TV+\frac12\delta^6\partial_X^3W
-\delta^7\partial_X(V^2)\\
&-2\delta\partial_X\bigl(e^{\delta^3W}-1-\delta^3W\bigr)
+\frac14\delta^9\partial_X\bigl((\partial_XW)^2\bigr).
\end{align*}
Expanding the exponential and collecting powers of $ \delta $ gives
\begin{align*}
Res_w={}&\delta^6\bigl(-\partial_TW_0-2\partial_XW_2\bigr)\\
&+\delta^7\bigl(-\partial_TW_1-2\partial_XW_3
-2W_0\partial_XW_0\bigr)
-2\delta^8\partial_X(W_0W_1)+ \mathcal{O}(\delta^9),\\
Res_v={}&\delta^6\bigl(-\partial_TW_0+2\partial_XW_2
+\tfrac12\partial_X^3W_0\bigr)\\
&+\delta^7\bigl(-\partial_TW_1+2\partial_XW_3
+\tfrac12\partial_X^3W_1-4W_0\partial_XW_0\bigr)\\
&+\delta^8\bigl(-\partial_TW_2-4\partial_X(W_0W_1)\bigr)
+ \mathcal{O}(\delta^9).
\end{align*}
At order $\delta^6$, the two equations require
\[
\partial_TW_0=-2\partial_XW_2,
\qquad
\partial_TW_0=2\partial_XW_2+\frac12\partial_X^3W_0.
\]
Adding them and then solving for $W_2$ yields
\begin{equation}\label{airyTruncatedLeading}
\partial_TW_0=\frac14\partial_X^3W_0,
\qquad W_2=-\frac18\partial_X^2W_0.
\end{equation}
At order $\delta^7$, we similarly obtain
\begin{align*}
\partial_TW_1&=-2\partial_XW_3-2W_0\partial_XW_0,\\
\partial_TW_1&=2\partial_XW_3+\frac12\partial_X^3W_1
-4W_0\partial_XW_0.
\end{align*}
Consequently,
\begin{equation}\label{airyTruncatedCorrection}
\partial_TW_1=\frac14\partial_X^3W_1
-\frac32\partial_X(W_0^2),
\qquad
W_3=-\frac18\partial_X^2W_1+\frac14W_0^2.
\end{equation}
The integration constants in $W_2,W_3$ are zero for 
Sobolev spaces on $\mathbb R$. We may prescribe $W_0(\cdot,0)$
and choose $W_1(\cdot,0)=0$.
The equation for $ W_1 $ can be solved with the variation of constants formula.
More explicitly, with $S(T)=\exp(T\partial_X^3/4)$,
\[
W_0(T)=S(T)W_0(0),\qquad
W_1(T)=-\frac32\int_0^T S(T-\sigma)
\partial_X\bigl(W_0(\sigma)^2\bigr)\,d\sigma.
\]
The coefficients through $\delta^7$ vanish. The residuals
are formally of order $\mathcal O(\delta^8)$. The spatial rescaling gives
$\mathcal O(\delta^{15/2})$ in the stated Sobolev norms.

For general $ \alpha > 2 $
the first corrections arise from nonlinearity and higher-order
dispersion. Their physical residuals have orders
$\delta^{2\alpha+1}$ and $\delta^{\alpha+5}$, respectively.
On the time scale $t=\mathcal O(\delta^{-3})$, these yield
corrections of orders $\delta^{2\alpha-2}$ and
$\delta^{\alpha+2}$. Thus, in the slow variables $(X,T)$,
the correction is of order
$\mathcal O(\delta^{\min\{2\alpha-2,\alpha+2\}})$.
Since $\alpha>2$, this is $o(\delta^\alpha)$.

\section{Higher order Burgers approximation}
\label{appHigherBurgers}

This section contains a proof of Lemma
\ref{lemresburgers} for the case $\alpha=1$. For this value we 
exemplarily construct a higher order Burgers
approximation
 with residual $\mathcal O(\delta^{9/2})
=\mathcal O(\delta^{2 \alpha+5/2})$ in $H^m\times H^{m+1}$.
We make the ansatz
\[
w_{\rm app}(x,t)=\delta W(X,T),\qquad
v_{\rm app}(x,t)=\delta V(X,T),
\]
with $ X=\delta(x-2t) $ and  $ T=\delta^2 t $. 
Substitution into \eqref{made1}-\eqref{made2} gives
\begin{align*}
Res_w&=2\delta^2\partial_X(W-V)-\delta^3\partial_TW
-2\delta^3V\partial_XW,\\
Res_v&=2\delta^2\partial_X(V-W)-\delta^3\partial_TV
+\frac12\delta^4\partial_X^3W-\delta^3\partial_X(V^2)\\
&\quad-2\delta\partial_X(e^{\delta W}-1-\delta W)
+\frac14\delta^5\partial_X((\partial_XW)^2).
\end{align*}
We set 
\[
W=\sum_{n=0}^2\delta^{n}W_n,
\qquad V=\sum_{n=0}^2\delta^{n}V_n,
\]
where the  $ V_j $ and $W_j$ are independent of $0 < \delta  \ll 1$.
Then the  $\mathcal{O}(\delta^2)$ terms imply $V_0=W_0$ if we  set the integration constant to zero. 

The determining equations for  $W_0,V_1$ 
can be found at $ \mathcal{O}(\delta^3) $. We find
\[
\partial_TW_0=-\frac32 \partial_X(W_0^2),
\qquad V_1=W_1+\frac14W_0^2.
\]
The first corrections $W_1,V_2$ can be found at $ \mathcal{O}(\delta^4) $.
We find
\[
\partial_TW_1+3\partial_X(W_0W_1)
=\frac14\partial_X^3W_0-\frac14\partial_X(W_0^3),
\]
\[
V_2=W_2+\frac12W_0W_1-\frac18\partial_X^2W_0+\frac1{24}W_0^3.
\]
Therefore, in the residuals terms of $ \mathcal{O}(\delta^5) $ remain.
The factor $\delta^{-1/2}$ comes from the scaling of the $ L^2 $-norm.

For  general $ \alpha \in (0,2) $ the first corrections arise from higher-order nonlinear terms
and dispersion. Their physical residuals have orders
$\delta^{3\alpha+1}$ and $\delta^{\alpha+3}$, respectively.
On the time scale $t=\mathcal O(\delta^{-(1+\alpha)})$,
these yield corrections of orders $\delta^{2\alpha}$ and
$\delta^2$. Thus, in the slow variables $(X,T)$,
the correction is of order
$\mathcal O(\delta^{\min\{2\alpha,2\}})$.
Since $0<\alpha<2$, this is $o(\delta^\alpha)$.


\medskip

{\bf Acknowledgements.}
The authors are grateful for discussions with Sarah Hofbauer, Xian Liao, and Robert Wegner.
The paper 
is partially supported by the Deutsche Forschungsgemeinschaft DFG through the SFB 1173 ''Wave phenomena'' Project-ID 258734477.

\bibliographystyle{alpha}
\bibliography{refs}

\end{document}